\documentclass[11pt]{amsart}
\usepackage{amssymb,amsmath,amsthm}
\newcommand{\leg}[2]{\genfrac{(}{)}{}{}{#1}{#2}}

\newtheorem{theorem}{Theorem}
\newtheorem{lemma}[theorem]{Lemma}

\newtheorem{proposition}[theorem]{Proposition}

\theoremstyle{remark}
\newtheorem*{remark}{Remark}
\theoremstyle{definition}
\newtheorem*{conjecture}{Conjecture}

\newtheorem*{example}{Example}

\numberwithin{theorem}{section} \numberwithin{equation}{section}
\newcommand{\ch}{{\text {\rm ch}}}

\newcommand{\C}{\mathbb{C}}

\newcommand{\Z}{\mathbb{Z}}

\newcommand{\SL}{{\text {\rm SL}}}

\newcommand{\Fp}{\mathbb{F}_p}

\newcommand{\calF}{\mathcal{F}}

\newcommand{\calW}{\mathcal{W}}

\begin{document}
\title[Wronskians of Andrews-Gordon series]
{Number theoretic properties of Wronskians of Andrews-Gordon
series}

\author{Antun Milas, Eric Mortenson, and Ken Ono}
\begin{abstract}
For positive integers $1\leq i\leq k$, we consider the arithmetic
properties of quotients of Wronskians in certain normalizations of
the Andrews-Gordon $q$-series
$$
\prod_{1\leq n\not \equiv 0,\pm i\pmod{2k+1}}\frac{1}{1-q^n}.
$$
This study is motivated by their appearance in conformal field
theory, where these series are essentially the irreducible
characters of $(2,2k+1)$ Virasoro minimal models. We  determine
the vanishing of such Wronskians, a result whose proof reveals
many partition identities. For example, if $P_{b}(a;n)$ denotes
the number of partitions of $n$ into parts which are not congruent
to $0, \pm a\pmod b$, then for every positive integer $n$ we have
$$
P_{27}(12; n)=P_{27}(6;n-1) + P_{27}(3;n-2).
$$
We also show that these quotients classify supersingular elliptic
curves in characteristic $p$. More precisely, if $2k+1=p$, where
$p\geq 5$ is prime, and the quotient is non-zero, then it is
essentially the locus of characteristic $p$ supersingular
$j$-invariants in characteristic $p$.
\end{abstract}

\thanks{The third author thanks the support of a
Packard Fellowship, a Romnes Fellowship, a Guggenheim
Fellowship, and grants from the National Science Foundation.}

\address{Department of Mathematics and Statistics,
SUNY Albany, Albany, New York 12222}
\email{amilas@math.albany.edu}

\address{Department of Mathematics, Penn State University,
University Park, Pennsylvania 16802}
\email{mort@math.psu.edu}

\address{Department of Mathematics, University of Wisconsin,
Madison, Wisconsin 53706}
\email{ono@math.wisc.edu}

 \maketitle

\section{Introduction and Statement of Results}

In two--dimensional conformal field theory and vertex operator
algebra theory (see \cite{Bo} \cite{FLM}), modular functions and modular
forms appear as graded dimensions, or characters, of infinite
dimensional irreducible modules. As a celebrated example, the graded dimension
of the Moonshine Module is the modular function
$$j(z)-744=q^{-1}+196884q+21493760q^2+\cdots
$$
(see
\cite{FLM}), where $q:=e^{2\pi i z}$ throughout.
Although individual characters are not always modular in this way, it can be
the case
that the vector spaces spanned by all of the irreducible characters
of a module are
invariant under the modular group \cite{Zhu}.
For example,
to construct an
automorphic form from an $\SL_2(\mathbb{Z})$--module,
one may simply
take the Wronskian of a basis of the module. In the case of the
Virasoro vertex operator algebras, such Wronskians were studied
by the first author (see
\cite{Milas1}, \cite{Milas2}) who obtained several classical $q$--series
identities related to modular forms using methods from representation theory.

Wronskian determinants in modular forms already play many
roles in number theory.
For example,
Rankin classified multi--linear
differential operators mapping automorphic forms to automorphic
forms
using Wronskians \cite{R}.
As another example,
the zeros of the Wronskian of a basis of weight
two cusp forms for a congruence subgroup $\Gamma_0(N)$ typically are
the Weierstrass points of the modular curve $X_0(N)$
\cite{FaKr}.

In view of this connection between Weierstrass points on modular
curves and Wronskians of weight 2 cusp forms, it is natural to
investigate the number theoretic properties of Wronskians of
irreducible characters. Here we make a first step in this
direction, and we consider an important class of models in vertex
operator algebra theory, those associated to $(2,2k+1)$ Virasoro
minimal models. These representations are important in conformal
field theory and mathematical physics and have been studied
extensively in the literature (see \cite{FF}, \cite{KW},
\cite{Milas2}, \cite{RC} and references therein).

We shall need some notation. Throughout, suppose that $k\geq
2$ is an integer. Define the rational number $c_{k}$ by
\begin{equation}\label{ck}
c_{k}:=1-\frac{3(2k-1)^2}{(2k+1)},
\end{equation}
and for each $1\leq i\leq k$, define $h_{i,k}$ by
\begin{equation}\label{hik}
h_{i,k}:=\frac{(2(k-i)+1)^2-(2k-1)^2}{8(2k+1)}.
\end{equation}
Let $L(c_k,h_{i,k})$ denote the irreducible lowest weight module
for the Virasoro algebra of central charge $c_k$ and weight
$h_{i,k}$ (see \cite{FF}, \cite{KW}, \cite{Milas1}). These
representations are $\mathbb{N}$--gradable and have finite
dimensional graded subspaces. Thus, we can define the formal
$q$--series
$${\rm dim}_{L(c_k,h_{i,k})}(q):=\sum_{n=0}^\infty {\rm dim}(L(c_k,h_{i,k})_n) q^n.$$
It is important to multiply the right hand side by the factor
$q^{h_{i,k}-\frac{c_k}{24}}$. The corresponding expression is
called the character of $L(c_k,h_{i,k})$ and will be denoted by
${\rm ch}_{i,k}(q)$. It turns out that (see \cite{RC}, \cite{KW})
\begin{equation}\label{ch}
\ch_{i,k}(q)=q^{(h_{i,k}-\frac{c_{k}}{24})}\cdot \prod_{1\leq n
\not \equiv 0, \pm i\pmod{2k+1}} \frac{1}{1-q^n}.
\end{equation}
Apart from the fractional powers of $q$ appearing in their
definition, such series have been studied extensively by Andrews,
Gordon, and of course Rogers and Ramanujan (for example, see \cite{A}). Indeed,
when $k=2$ we have that $c_{5}=-\frac{22}{5}$, and the two
corresponding characters are essentially the products appearing in
the celebrated Rogers-Ramanujan identities
\begin{displaymath}
\begin{split}
\ch_{1,2}(q)&=q^{\frac{11}{60}}
\prod_{n\geq 0}\frac{1}{(1-q^{5n+2})(1-q^{5n+3})},\\
\ch_{2,2}(q)&=q^{-\frac{1}{60}}
\prod_{n\geq 0}\frac{1}{(1-q^{5n+1})(1-q^{5n+4})}.
\end{split}
\end{displaymath}

In view of the discussion above, we
investigate, for each $k$, the
Wronskians for the complete sets
of characters
$$\{\ch_{1,k}(q), \ch_{2,k}(q),\dots, \ch_{k,k}(q)\}.
$$
For each $k\geq 2$, define $\calW_k(q)$ and
$\calW'_k(q)$ by
\begin{equation}\label{W}
\calW_k(q):=\alpha(k)\cdot \det \left( \begin{matrix}
\ch_{1,k} & \ch_{2,k} &\cdots & \ch_{k,k}\\
\ch_{1,k}' & \ch_{2,k}' &\cdots & \ch_{k,k}'\\
\vdots     &  \vdots    & \vdots & \vdots \\
\ch_{1,k}^{(k-1)} & \ch_{2,k}^{(k-1)} & \cdots & \ch_{k,k}^{(k-1)}
\end{matrix}
\right),
\end{equation}
\begin{equation}\label{Wp}
\calW_k'(q):=\beta(k)\cdot \det \left( \begin{matrix}
\ch_{1,k}' & \ch_{2,k}' &\cdots & \ch_{k,k}'\\
\ch_{1,k}^{(2)} & \ch_{2,k}^{(2)} &\cdots & \ch_{k,k}^{(2)}\\
\vdots     &  \vdots    & \vdots & \vdots \\
\ch_{1,k}^{(k)} & \ch_{2,k}^{(k)} & \cdots & \ch_{k,k}^{(k)}
\end{matrix}
\right).
\end{equation}
Here $\alpha(k)$ (resp. $\beta(k)$) is chosen so that the
$q$-expansion of $\calW_k(q)$ (resp. $\calW_k^{\prime}(q)$) has
leading coefficient $1$ (resp. $1$ or $0$), and differentiation
is given by
$$
\left(\sum a(n)q^n\right)' :=\sum na(n)q^n,
$$
which equals $\frac{1}{2\pi i}\cdot \frac{d}{dz}$ when $q:=e^{2\pi
i z}$.

It turns out that $\calW_k(q)$ is
easily described in terms of Dedekind's eta-function
(see also Theorem~6.1 of \cite{Milas2}),
which is defined for $z\in \mathbb{H}$, $\mathbb{H}$ denoting the usual upper half-plane of $\mathbb{C}$, by
$$
\eta(z):=q^{\tfrac{1}{24}}\prod_{n=1}^{\infty}(1-q^n).
$$
\begin{theorem}\label{denwronsk}
If $k\ge2$, then
\begin{displaymath}
\calW_k(q)=\eta(z)^{2k(k-1)}.
\end{displaymath}
\end{theorem}
Instead of directly computing $\calW_k^{\prime}(q)$, we
investigate the quotient
\begin{equation}\label{Fk}
\calF_k(z):=\frac{\calW_k'(q)}{\calW_k(q)}.
\end{equation}
It turns out that these $q$-series $\calF_k(z)$ are modular forms
of weight $2k$ for the full modular group $\SL_2(\Z)$. To make
this more precise, suppose that $E_4(z)$ and $E_6(z)$ are the
standard Eisenstein series
\begin{equation}
E_4(z)=1+240\sum_{n=1}^{\infty}\sum_{d\mid n}d^{3} q^n\ \ \ \ \
{\text {\rm and}}\ \ \ \ \
E_6(z)=1-504\sum_{n=1}^{\infty}\sum_{d\mid n}d^{5} q^n,
\end{equation}
and that $\Delta(z)$ and $j(z)$ (as before) are the usual modular forms
\begin{equation}
\Delta(z)=\frac{E_4(z)^3-E_6(z)^2}{1728}\ \ \ \ \ {\text {\rm and}}\ \ \ \ \
j(z)=\frac{E_4(z)^3}{\Delta(z)}.
\end{equation}
\begin{theorem}\label{milastheorem}
If $k\geq 2$, then $\calF_k(z)$ is a weight $2k$ holomorphic modular form
on $\SL_2(\Z)$.
\end{theorem}

\begin{example}
For $2\leq k\leq5$, it turns out that
\begin{displaymath}
\begin{split}
\calF_2(z)&=E_4(z), \ \ \  \
\ \ \calF_3(z)=E_6(z),\\
\calF_4(z)&=E_4(z)^2,\ \ \ \ \
\calF_5(z)=E_4(z)E_6(z).
\end{split}
\end{displaymath}
\end{example}

Some of these modular forms are identically
zero. For example, we have that $\calF_{13}(z)=0$,
a consequence of the $q$-series identity
\begin{equation}\label{k13}
\ch_{12,13}(q)-\ch_{6,13}(q)-\ch_{3,13}(q)=1,
\end{equation}
which will be proved later.  The following result completely determines those $k$ for which
$\calF_k(z)=0$.

\begin{theorem}\label{vanwronsk}
The modular form
$\calF_k(z)$ is identically zero precisely for those
$k$ of the form $k=6t^2-6t+1$ with $t\geq 2$.
\end{theorem}

Since the Andrews-Gordon series are the partition generating functions
$$
\sum_{n=0}^{\infty}P_b(a;n)q^n=
\prod_{1\leq n\not \equiv \pm 0, a\pmod{b}}\frac{1}{1-q^n},
$$
$q$-series identity (\ref{k13}) implies, for positive $n$,
the shifted partition identity
$$
P_{27}(12; n)=P_{27}(6;n-1) + P_{27}(3;n-2).
$$
This partition identity is a special case of the following
theorem which is a corollary to the proof of
Theorem~\ref{vanwronsk}.

\begin{theorem}\label{partidentity}
If $t\geq 2$, then for every positive integer $n$ we have
\begin{displaymath}
P_{b(t)}(a^{-}(t,0);n)
=\sum_{r=1}^{t-1}(-1)^{r+1}
\left( P_{b(t)}(a^{+}(t,r);n-\omega^{-}(r))
+P_{b(t)}(a^{-}(t,r);
n-\omega^{+}(r))\right),
\end{displaymath}
where
\begin{displaymath}
\begin{split}
a^{-}(t,r)&:=(2t-1)(3t-3r-2),\\
a^{+}(t,r)&:=(2t-1)(3t-3r-1),\\
b(t)&:=3(2t-1)^2,\\
\omega^{-}(r)&:=(3r^2-r)/2,\\
\omega^{+}(r)&:=(3r^2+r)/2.
\end{split}
\end{displaymath}
\end{theorem}

The forms $\calF_k(z)$ also provide deeper number theoretic
information.
Some of them parameterize isomorphism classes
of supersingular elliptic curves in characteristic $p$. To make this
precise, suppose that $K$ is a field with characteristic $p>0$, and
let $\overline{K}$ be its algebraic closure.
An elliptic curve $E$ over $K$ is {\it supersingular} if the group
$E(\overline{K})$ has no $p$-torsion.
This connection is phrased in terms of ``divisor polynomials"
of modular forms.

We now describe these polynomials.
If $k\ge4$ is even, then define $\tilde E_k(z)$ by
\begin{equation}\label{Ektilde}
\tilde E_k(x):= \begin{cases}
1\ \ \ \  &{\text {\rm if}}\ k\equiv 0 \pmod{12},\\
E_4(z)^2E_6(z)\ \ \ \  &{\text {\rm if}}\ k\equiv 2 \pmod{12},\\
E_4(z)\ \ \ \ &{\text {\rm if}}\ k\equiv 4 \pmod{12},\\
E_6(z)\ \ \ \ &{\text {\rm if}}\ k\equiv 6 \pmod{12},\\
E_4(z)^2\ \ \ \ &{\text {\rm if}}\ k\equiv 8 \pmod{12},\\
E_4(z)E_6(z)\ \ \ \ &{\text {\rm if}}\ k\equiv 10 \pmod{12}
\end{cases}
\end{equation}
(see Section $2.6$ of \cite{cbms} for further details on divisor polynomials).
As usual, let $M_k$ denote the space of holomorphic
weight $k$ modular forms on $\SL_2(\Z)$.
If we write $k$ as
\begin{equation}
k=12m+s \ \ \text{with }s\in\{0,4,6,8, 10,14\},
\end{equation}
then $\dim_{\C}(M_k)=m+1$,  and every modular form $f(z)\in M_k$
factorizes as
\begin{equation}
f(z)=\Delta(z)^m\tilde E_k(z)\tilde F(f,j(z)),
\end{equation}
where $\tilde F$ is a polynomial of degree $\le m$ in $j(z)$.
Now define the polynomial $h_k(x)$ by
\begin{equation}\label{hk}
h_k(x):= \begin{cases}
1\ \ \ \  &{\text {\rm if}}\ k\equiv 0 \pmod{12},\\
x^2(x-1728)\ \ \ \ &{\text {\rm if}}\ k\equiv 2 \pmod{12},\\
x\ \ \ \ &{\text {\rm if}}\ k\equiv 4 \pmod{12},\\
x-1728\ \ \ \ &{\text {\rm if}}\ k\equiv 6 \pmod{12},\\
x^2\ \ \ \ &{\text {\rm if}}\ k\equiv 8 \pmod{12},\\
x(x-1728)\ \ \ \ &{\text {\rm if}}\ k\equiv 10 \pmod{12}.
\end{cases}
\end{equation}
If $f(z)\in M_k$, then define the {\it{divisor polynomial}} $F(f,x)$ by
\begin{equation}\label{dp}
F(f,x):=h_k(x)\tilde F(f,x).
\end{equation}
If $j(E)$ denotes the usual $j$-invariant
of an elliptic curve $E$, then the
{\it characteristic $p$ locus of supersingular
$j$-invariants} is the polynomial in $\mathbb{F}_p[x]$ defined by
\begin{equation}\label{slocus}
S_p(x):= \prod_{E/\overline{\mathbb{F}}_p\ {\text {\rm
supersingular}}} (x-j(E)),
\end{equation}
the product being over isomorphism classes of supersingular
elliptic curves.
The following congruence modulo $37$, when $k=18$, is a special
case of our general result
\begin{displaymath}
\begin{split}
F(\calF_{18},j(z))&=j(z)^3-\frac{2^{13}\cdot 3^4\cdot 89\cdot
1915051410991641479}{17\cdot 43\cdot 83\cdot 103\cdot 113\cdot 163\cdot
523 \cdot 643 \cdot 919\cdot 1423}\cdot j(z)^2+\cdots\\
&\equiv (j(z)+29)(j(z)^2+31j(z)+31) \pmod{37}\\
&= S_{37}(j(z)).
\end{split}
\end{displaymath}
The following result provides a general class of $k$ for which
$F(\calF_k,j(z))\pmod p$ is the supersingular locus $S_p(j(z))$.

\begin{theorem}\label{sslcong}
If $p\geq 5$ is prime, and
$k=(p-1)/2$ is not of the form $6t^2-6t+1$, where $t\ge2$, then
\begin{displaymath}
F(\calF_k,j(z))\equiv  S_p(j(z)) \pmod p.
\end{displaymath}
\end{theorem}

\medskip

In Section~\ref{sect:6} we prove Theorem
~\ref{milastheorem}.  In Section~\ref{sect:1} we provide
the preliminaries required for the proofs of
Theorems~\ref{vanwronsk} and ~\ref{sslcong}.  These theorems,
along with Theorem~\ref{partidentity},
are then proved in
Sections~\ref{sect:2} and ~\ref{sect:3} respectively.
In Section~\ref{sect:7} we give a conjecture
concerning the zeros of $\widetilde F(\calF_k,x)$.

\medskip


\section{Differential operators and the proof of Theorem~\ref{milastheorem}}
\label{theorem1}\label{sect:6}

Here we study the quotient
\begin{equation}\label{FFK}
\calF_k(q):=\frac{\mathcal{W}'_k(q)}{\mathcal{W}_k(q)},
\end{equation}
and prove Theorem~\ref{milastheorem}.
Needless to say, the previous expression can be defined for an
arbitrary set  $\{f_1(q),...,f_k(q)\}$ of holomorphic functions in
$\mathbb{H}$, where each $f_i(q)$ has a $q$--expansion. In the
$k=1$ case, (\ref{FFK}) is just the logarithmic derivative of
$f_1$. We also introduce generalized Wronskian determinants. For
$0 \leq i_1 < i_2 < \cdots < i_k,$ let
\begin{equation}\label{genwronsk}
{W}^{i_1,...,i_k}(f_1,...,f_k)(q)={\rm det} \left( \begin{array}{ccccc} f^{(i_1)}_1 & f^{(i_1)}_2 & . & . & f^{(i_1)}_{k} \\
f^{(i_2)}_1 & f^{(i_2)}_2 & . & . & f^{(i_2)}_{k} \\ . & . & . & . & . \\
. & . & . & . & . \\ f_1^{(i_k)} & f_2^{(i_k)} & . & . &
f_{k}^{(i_k)} \end{array} \right).
\end{equation}
Clearly,
${W}^{0,1,...,k-1}(f_1,...,f_k)$ is the ordinary Wronskian. We
will write $\mathcal{W}^{i_1,...,i_k}(f_1,...,f_k)$ for the
normalization of $W^{i_1,...,i_k}(f_1,...,f_k) \neq 0$, where the
leading coefficient in the $q$-expansion is one.

Before we prove Theorem~\ref{milastheorem}, we recall
the definition of the ``quasi-modular" form
\begin{equation}\label{E2}
E_2(z)=1-24 \sum_{n=1}^\infty \sum_{d\mid n} d q^n,
\end{equation}
which plays an important role in the proof of
the following result (also see \cite{Milas1}).
\begin{theorem} \label{Milas0}
For $k \geq 2$, the set $\{\ch_{1,k}(q),...,\ch_{k,k}(q) \}$ is a
fundamental system of a $k$-th order linear differential equation
of the form
\begin{equation} \label{system}
\left(q \frac{d}{dq} \right)^k y+\sum_{i=0}^{k-1} P_i(q) \left( q
\frac{d}{dq} \right)^i y=0,
\end{equation}
where $P_i(q) \in \mathbb{Q}[E_2,E_4,E_6]$. Moreover, $P_0(q)$ is
a modular form for $\SL_2(\mathbb{Z})$.
\end{theorem}
\begin{proof}
The first part was already proven in \cite{Milas2}, Theorem 6.1.
Let
\begin{equation} \label{thetafor}
\Theta_k:=\left(q \frac{d}{dq}\right)-\frac{k}{12} E_2(z)
\end{equation}
It is well known that $\Theta_k$ sends a modular form of weight
$k$ to a modular form of weight $k+2$. Then Theorem 5.3 in
\cite{Milas1} implies that the equation (\ref{system}) can be
rewritten as
$$\Theta^k y + \sum_{i=1}^{k-1}  Q_i(q) \Theta^i y+Q_0(q)y =0,$$ where
for $i \geq 1$
$$\Theta^i:=\Theta_{2i-2} \circ \cdots \circ
\Theta_2 \circ \Theta_0,$$ and
$$Q_i(q) \in \mathbb{Q}[E_4,E_6].$$
Clearly, $P_0(q)=Q_0(q)$. The proof follows.
\end{proof}

\begin{proof}[Proof of Theorem ~\ref{milastheorem}]
Let $\{f_1(q),...,f_k(q) \}$ be a linearly independent set of
holomorphic functions in the upper half-plane. Then there is a
unique $k-$th
 order linear differential operator with meromorphic
coefficients
$$P=\left(q \frac{d}{dq} \right)^k+\sum_{i=0}^{k-1}P_i(q)\left( q
\frac{d}{dq}\right)^i,$$ such that $\{f_1(q),...,f_k(q)\}$ is a
fundamental system of
$$P(y)=0.$$
Explicitly,
\begin{equation} \label{correspond}
P(y)=(-1)^{k}
\frac{W^{0,1,...,k}(y,f_1,...,f_k)}{W^{0,1,...,k-1}(f_1,...,f_k)}.
\end{equation}
In particular,
$$P_0(q)=(-1)^k \frac{{W}^{1,2,...,k}(f_1,...,f_k)}{{W}^{0,1,...,k-1}(f_1,...,f_k)}.$$
Thus,
$$P_0(q)= \lambda_k \frac{\mathcal{W}'_k(q)}{\mathcal{W}_k(q)},$$
for some nonzero constant $\lambda_k$. Now, we specialize
$f_i(q)={\rm ch}_{i,k}(q)$ and apply Theorem \ref{Milas0}.
\end{proof}

\begin{remark}
The techniques from \cite{Milas1} can be used to give explicit
formulas for $\frac{\mathcal{W}'_k(q)}{\mathcal{W}_k(q)}$ in terms
of Eisenstein series. However, this computation becomes very
tedious for large $k$.
\end{remark}

\section{Preliminaries for Proofs of Theorems~\ref{vanwronsk} and~\ref{sslcong}}\label{sect:1}

In this section we recall essential preliminaries regarding $q$-series
and divisor polynomials of modular forms.

\subsection{Classical $q$-series identities}

We begin by recalling Jacobi's triple product identity and Euler's
pentagonal number theorem.

\begin{theorem}\label{jtpid}
(Jacobi's Triple Product Identity) For $y\ne 0$ and $|q|<1$, we
have
$$
\sum_{n=-\infty}^{\infty}y^nq^{n^2}=\prod_{n=1}^{\infty}(1-q^{2n})(1+yq^{2n-1})(1+y^{-1}q^{2n-1}).
$$
\end{theorem}

\begin{theorem}\label{epnt}
(Euler's Pentagonal Number Theorem)
The following $q$-series identity is true:
$$
\prod_{n=1}^{\infty}(1-q^n)=\sum_{m=-\infty}^{\infty}(-1)^m q^{\tfrac12m(3m-1)}
$$
\end{theorem}

\subsection{Divisor polynomials and Deligne's theorem}

If $p\geq 5$ is prime, then the
supersingular loci $S_p(x)$ and $\widetilde{S}_p(x)$ are defined in $\Fp[x]$
by the following products over isomorphism classes of supersingular
elliptic curves:
$$
S_p(x):= \prod_{E/\overline{\mathbb{F}}_p\ {\text {\rm
supersingular}}} (x-j(E)),
$$
\begin{equation}\label{ssplocus}
\widetilde{S}_p(x):= \prod_{\substack{E/\overline{\mathbb{F}}_p \
\ {\text {\rm supersingular}}\\ j(E) \not \in \{0, 1728\}}}
(x-j(E)).
\end{equation}

For such primes $p$,
let $\mathfrak{S}_p$ denote the set of those supersingular
$j$-invariants in characteristic $p$ which are in $\Fp-\{0,
1728\}$, and  let $\mathfrak{M}_p$ denote the set of monic
irreducible quadratic polynomials in $\Fp[x]$ whose roots are
supersingular $j$-invariants.
The polynomial $S_p(x)$ splits completely in
$\mathbb{F}_{p^2}$ (\cite{Silverman}).
Define $\epsilon_{\omega}(p)$ and $\epsilon_{i}(p)$ by
$$
\epsilon_{\omega}(p):=\begin{cases}
0 &{\text {\rm if}}\ p\equiv 1 \pmod 3,\\
1 &{\text {\rm if}}\ p\equiv 2 \pmod 3,
\end{cases}
$$
$$
\epsilon_{i}(p):=\begin{cases}
0 &{\text {\rm if}}\ p\equiv 1 \pmod 4,\\
1 &{\text {\rm if}}\ p\equiv 3 \pmod 4,
\end{cases}
$$
The following proposition relates
$S_p(x)$ to $\tilde S_p(x)$ (\cite{Silverman}).

\begin{proposition}
If $p\ge5$ is prime, then
\begin{displaymath}
\begin{split}
S_p(x)&=x^{\epsilon_{\omega}(p)}(x-1728)^{\epsilon_{i}(p)}\cdot\prod_{\alpha\in \mathfrak{S}_p}(x-\alpha)\cdot\prod_{g\in \mathfrak{M}_p}g(x)\\
&=x^{\epsilon_{\omega}(p)}(x-1728)^{\epsilon_{i}(p)}\tilde S_p(x).
\end{split}
\end{displaymath}
\end{proposition}

Deligne found the following explicit description of these polynomials (see \cite{dwork}, \cite{Serre}).
\begin{theorem}\label{deltheorem}
If $p\ge5$ is prime, then
$$
F(E_{p-1},x)\equiv S_p(x) \pmod p.
$$
\end{theorem}

\begin{remark} In a beautiful paper \cite{KanekoZagier}, Kaneko and Zagier
provide a simple proof of Theorem~\ref{deltheorem}.
\end{remark}

\begin{remark}
The Von-Staudt congruences imply for primes $p$, that $\tfrac{2(p-1)}{B_{p-1}}
\equiv 0 \pmod p$,
where $B_n$ denotes the usual $n$th Bernoulli number.  It follows that if
$$E_k(z)=1-\frac{2k}{B_k}\sum_{n=1}^{\infty}\sum_{d\mid n}d^{k-1}q^n
$$
is the usual weight $k$ Eisenstein series, then
\begin{displaymath}
E_{p-1}(z)\equiv 1 \pmod p.
\end{displaymath}
If $p\geq 5$ is prime, then Theorem~\ref{deltheorem}
combined with the definition of divisor
polynomials, implies
that if $f(z)\in M_{p-1}$ and $f(z)\equiv 1 \pmod p$, then
\begin{displaymath}
F(f,j(z))\equiv S_p(j(z)) \pmod p.
\end{displaymath}
\end{remark}

\section{The Vanishing of $\calF_k(z)$}\label{sect:2}

Here we prove
Theorem~\ref{vanwronsk}
and Theorem~\ref{partidentity}.
To prove these results, we first
require some notation and two  technical lemmas.
For simplicity, we will write $\ch_{i}(q)$
for $\ch_{i,k}(q)$ when $k$ is understood.

We define
\begin{equation}
\Theta(y,q):=\sum_{n=-\infty}^{\infty}y^nq^{n^2},
\end{equation}
and we consider the sum
\begin{equation}
A_t(q):=\Theta(-q^{\tfrac12(2t-1)},q^{\tfrac32(2t-1)^2})+\sum_{r=1}^{t-1}(-1)^r\Psi^{-}_{r,t}(q)+\sum_{r=1}^{t-1}(-1)^r\Psi^{+}_{r,t}(q),
\end{equation}
where
$$
\Psi^{-}_{r,t}(q):=q^{\tfrac12r(3r-1)}\Theta(-q^{\tfrac12(6r-1)(2t-1)},q^{\tfrac32(2t-1)^2}),
$$
and
$$
\Psi^{+}_{r,t}(q):=q^{\tfrac12r(3r+1)}\Theta(-q^{\tfrac12(6r+1)(2t-1)},q^{\tfrac32(2t-1)^2}).
$$

\begin{lemma}\label{lem1}
If $t\ge2$ and $k=6t^2-6t+1$, then we have the following
$q$-series identity
\begin{displaymath}
\begin{split}
\prod_{n=1}^{\infty}(1-q^n)^{-1}\cdot A_t(q)=&
\ch_{((2t-1)(3t-2))}(q)\\
&+\sum_{r=1}^{t-1}(-1)^r\ch_{((2t-1)(3t-3r-1))}(q)
+\sum_{r=1}^{t-1}(-1)^r\ch_{((2t-1)(3t-3r-2))}(q).
\end{split}
\end{displaymath}
\end{lemma}

\begin{proof}
We examine the summands  in $A_t(q)$.
Using Theorem~\ref{jtpid}, the first term is
\begin{displaymath}
\begin{split}
\Theta&(-q^{\tfrac12(2t-1)},q^{\tfrac32(2t-1)^2})\\
=&\prod_{n=1}^{\infty}(1-q^{3(2t-1)^2n})(1-q^{\tfrac12(2t-1)+\tfrac32(2t-1)^2(2n-1)})(1-q^{-\tfrac12(2t-1)+\tfrac32(2t-1)^2(2n-1)})\\
=&\prod_{n=1}^{\infty}(1-q^{3(2t-1)^2n})(1-q^{-(\tfrac32(2t-1)^2-\tfrac12(2t-1))+3(2t-1)^2n})(1-q^{\tfrac32(2t-1)^2-\tfrac12(2t-1)+3(2t-1)^2(n-1)}).
\end{split}
\end{displaymath}
Noting that $\frac12(3(2t-1)^2-(2t-1))=(2t-1)(3t-2)$,
we have that
\begin{displaymath}
\Theta(-q^{\tfrac12(2t-1)},q^{\tfrac32(2t-1)^2})
=\prod_{n=1}^{\infty}(1-q^n)\cdot \ch_{((2t-1)(3t-2))}(q).
\end{displaymath}
Arguing with  Theorem \ref{jtpid} again,
we find that
\begin{displaymath}
\Psi^{-}_{r,t}(q)=\prod_{n=1}^{\infty}(1-q^n)\cdot
\ch_{((2t-1)(3t-3r-1))}(q),
\end{displaymath}
and
\begin{displaymath}
\Psi^{+}_{r,t}(q)=\prod_{n=1}^{\infty}(1-q^n)\cdot
\ch_{((2t-1)(3t-3r-2))}(q).
\end{displaymath}
The lemma follows easily.
\end{proof}

\begin{lemma}\label{lem2}
If $t\ge2$, then we have the following $q$-series identity
$$
A_t(q)=\prod_{n=1}^{\infty}(1-q^n).
$$
\end{lemma}

\begin{proof}
It suffices to show that
$A_t(q)$ is the $q$-series in
Euler's Pentagonal Number Theorem.
Write $A_t(q)$
in a more recognizable form beginning with the first term in $A_t(q)$
\begin{displaymath}
\begin{split}
\Theta(-q^{\tfrac12(2t-1)},q^{\tfrac32(2t-1)^2})&=\sum_{n=-\infty}^{\infty}(-1)^nq^{\tfrac12(2t-1)n+\tfrac32(2t-1)^2n^2}\\
&=\sum_{n=-\infty}^{\infty}(-1)^nq^{\tfrac12(2t-1)n(3(2t-1)n-1)},
\end{split}
\end{displaymath}
where we have replaced $n$ by $-n$ in the final sum.

For $\Psi^{-}_{r,t}(q)$ and $\Psi^{+}_{r,t}(q)$ we have
\begin{displaymath}
\begin{split}
(-1)^r\Psi^{-}_{r,t}(q)&=\sum_{n=-\infty}^{\infty}(-1)^{n+r}q^{\tfrac12r(3r-1)+\tfrac12(6r-1)(2t-1)n+\tfrac32(2t-1)^2n^2}\\
&=\sum_{n=-\infty}^{\infty}(-1)^{(2t-1)n+r}
q^{\tfrac12((2t-1)n+r)(3((2t-1)n+r)-1)},
\end{split}
\end{displaymath}
and we have
\begin{displaymath}
\begin{split}
(-1)^r\Psi^{+}_{r,t}(q)&=\sum_{n=-\infty}^{\infty}(-1)^{n+r}q^{\tfrac12r(3r+1)+\tfrac12(6r+1)(2t-1)n+\tfrac32(2t-1)^2n^2}\\
&=\sum_{n=-\infty}^{\infty}(-1)^{n+r}q^{\tfrac12((2t-1)n+r)(3((2t-1)n+r)+1)}\\
&=\sum_{n=-\infty}^{\infty}(-1)^{(2t-1)n+(2t-1)-r}q^{\tfrac12 ((2t-1)n+(2t-1)-r) (3((2t-1)n+(2t-1)-r)-1) },
\end{split}
\end{displaymath}
where in the last line we substituted $-n$ for $n$ and then $n+1$ for $n$.
By combining these series, the claim
follows easily from Theorem \ref{epnt}.
\end{proof}

\begin{proof}[Proof of Theorem~\ref{vanwronsk}]
We first show that if $k=6t^2-6t+1$, then $\calW^{\prime}_k(q)$ vanishes.  Recalling that $2k+1=3(2t-1)^2$ and using the above lemmas, we obtain
\begin{displaymath}
{\ch}_{((2t-1)(3t-2))}(q)
+\sum_{r=1}^{t-1}(-1)^r\ch_{((2t-1)(3t-3r-1))}(q)
+\sum_{r=1}^{t-1}(-1)^r\ch_{((2t-1)(3t-3r-2))}(q)=1.
\end{displaymath}
This gives us a linear relationship among the columns, and the Wronskian is then identically zero.

Define, for $1\leq i\leq k$, the rational number
\begin{displaymath}
a(i,k):=h_{i,k}-\tfrac{c_{k}}{24}.
\end{displaymath}
If $k\ne 6t^2-6t+1$, it is straightforward to show that $a(i,k)\ne0$.
Noting that
\begin{displaymath}
\ch_{i,k}(q)=q^{a(i,k)}+\dots
\end{displaymath}
and making a few simple observations such as $a(i,k)>a(i+1,k)$ and $a(k,k)>-1$, it follows that $\calW^{\prime}_k(q)$ cannot vanish.
Specifically, if the Wronskian vanished, then
we would have a linear dependence of the characters.
However, because $a(i,k)\ne0$ and $a(i,k)>-1$, this is not possible.
\end{proof}

\begin{proof}[Proof of Theorem~\ref{partidentity}]
For $t\geq 2$ and
$k=6t^2-6t+1$, the proof of Theorem~\ref{vanwronsk} gives the
following identity
\begin{displaymath}
\begin{split}
\prod_{n\not\equiv 0, \pm (2t-1)(3t-2)\pmod {2k+1}}^{\infty}
&\frac1{1-q^n}+\sum_{r=1}^{t-1}(-1)^rq^{\tfrac12r(3r-1)}\cdot
\prod_{n\not\equiv 0, \pm (2t-1)(3t-3r-1)\pmod{2k+1}}^{\infty}
\frac{1}{1-q^n}\\
&+\sum_{r=1}^{t-1}(-1)^rq^{\tfrac12r(3r+1)}\cdot\prod_{n\not\equiv 0,
\pm (2t-1)(3t-3r-2) \pmod{2k+1}}^{\infty}\frac{1}{1-q^n}=1.
\end{split}
\end{displaymath}
The proof now follows by inspection.
\end{proof}

\section{Supersingular Polynomial Congruences}\label{sect:3}

Here we prove Theorem~\ref{sslcong}:  the congruence
\begin{displaymath}
F(\calF_k,j(z))\equiv  S_p(j(z)) \pmod p,
\end{displaymath}
which holds for primes $5\leq p=(2k+1)$, where
$k\neq 6t^2-6t+1$ with $t\geq 2$.

We begin with a technical lemma.
\begin{lemma}\label{wronsklem}
For $k$ with $2k+1=p$, $p$ a prime, and $k$ not of the form $6t^2-6t+1$, $t\ge2$, then
$\calF_k(z)$ has $p$-integral coefficients, and satisfies the congruence
\begin{displaymath}
\calF_k(z)\equiv 1 \pmod p.
\end{displaymath}
\end{lemma}
\begin{proof}
If we expand $\calW^{\prime}_k(q)$ by minors along its bottom row, and if we expand $\calW_k(q)$ by minors along its top row, we have that the quotient of the Wronskians is the normalization of
\begin{displaymath}
\frac{\displaystyle{\sum_{i=1}^{k}\ch_{i,k}^{(k)}(q)\det(M_i)}}{\displaystyle{\sum_{i=1}^{k}\ch_{i,k}(q)\det(M_i)}},
\end{displaymath}
where the $M(i)$'s are the respective minors.  We fix an $i$ and consider the term
\begin{displaymath}
\ch_{i,k}^{(k)}(q)=\sum_{n=0}^{\infty}(n+a(i,k))^kb_{i,k}(n)q^{n+a(i,k)},
\end{displaymath}
where
\begin{displaymath}
\ch_{i,k}(q)=\sum_{n=0}^{\infty}b_{i,k}(n)q^{n+a(i,k)}.
\end{displaymath}
We note that
\begin{displaymath}
a(i,k)=\frac{(2k+1)(3k+1-6i)+6i^2}{12(2k+1)}.
\end{displaymath}
Multiplying the numerator by $(12(2k+1))^k$ to clear out the denominators in the $a(i,k)$'s, the quotient of the Wronskians is then just the normalization of
\begin{displaymath}
\frac{\displaystyle{\sum_{i=1}^{k}\sum_{n=0}^{\infty}
\left(12(2k+1)n +(2k+1)(3k+1-6i)+6i^2 \right)^k
b_{i,k}(n)q^{n+a(i,k)}
\det(M_i)}}{\displaystyle{\sum_{i=1}^{k}\ch_{i,k}(q)\det(M_i)}}
\end{displaymath}
However if we compute this modulo $p$, and note that
$k=\tfrac{p-1}2$, we have
\begin{displaymath}
\left(12(2k+1)n+(2k+1)(3k+1-6i)+6i^2 \right)^k\equiv \leg{6}{p}
\pmod p.
\end{displaymath}
The $p$-integrality follows from Theorem~\ref{denwronsk}, and it then follows that modulo $p$,
the quotient is just $1$.
\end{proof}

\begin{proof}[Proof of Theorem~\ref{vanwronsk}]
Here we simply combine Theorem~\ref{milastheorem}, Lemma~\ref{wronsklem}
and the second remark at the end of Section~\ref{sect:1}.
\end{proof}
\begin{remark}
There are cases for which
\begin{displaymath}
\calF_k(z)\equiv 1 \pmod p,
\end{displaymath}
with $2k+1\neq p$.
It would interesting to completely determine all the conditions for which
such a congruence holds.
By the theory of modular forms `mod $p$',
it follows that such $k$ must have the property that
$2k=a(p-1)$ for some positive integer $a$.
A resolution of this problem requires determining conditions
for which $\calF_k(z)$ has $p$-integral coefficients, and also
the extra conditions guaranteeing the above congruence.
\end{remark}
\begin{remark}
The methods of this paper can be used to reveal many more congruences
relating supersingular $j$-invariants to the divisor
polynomials $F(\calF_k,j)\pmod p$.
For example, if
$$(k,p) \in \{(10,17), (16,29), (17,31), (22,41), (23,43),
(28,53)\},
$$
we have that
\begin{displaymath}
F(\calF_k,j(z))\equiv j(z)\cdot S_p(j(z)) \pmod p.
\end{displaymath}
Such congruences follow from the multiplicative structure
satisfied by divisor polynomials as described in
Section 2.8 of \cite{cbms}.
\end{remark}

\section{A conjecture on the zeros of $F(\calF_k,x)$}\label{sect:7}

Our Theorem~\ref{sslcong} shows that the divisor polynomial modulo
$p$, for certain $\calF_k(z)$, is the locus of supersingular
$j$-invariants in characteristic $p$. We proved this theorem by
showing that
$$
\calF_{\frac{p-1}{2}}(z)\equiv E_{p-1}(z)\equiv 1\pmod p,
$$
and we then obtained the desired conclusion by
applying a famous result of Deligne.
In view of such close relationships between certain
$\calF_k(z)$ and $E_{2k}(z)$, it is natural to investigate other
properties of $E_{2k}(z)$ which may be shared by
$\calF_k(z)$.
A classical result of Rankin and Swinnerton-Dyer proves that
every $F(E_{2k},x)$ has simple roots, all of which are
real and lie in the interval $[0,1728]$.
Numerical evidence strongly supports the following conjecture.

\begin{conjecture} If $k\geq 2$ is a positive integer for which
$k\ne 6t^2-6t+1$ with $t\ge2$, then
$F(\calF_k,x)$ has simple roots, all of which are
real and are in the interval $[0, 1728]$.
\end{conjecture}

\end{document}